# Robust Transceiver Optimization for Downlink Multiuser MIMO Systems

Tadilo Endeshaw, Batu Krishna Chalise and Luc Vandendorpe

*Abstract*— This paper addresses the joint transceiver design for downlink multiuser multiple-input multiple-output (MIMO) systems, with imperfect channel state information (CSI) at the base station (BS) and mobile stations (MSs). By incorporating antenna correlation at both ends of the channel and taking channel estimation errors into account, we solve two robust design problems: minimization of the weighted sum mean-square-error (MSE) and minimization of the maximum weighted MSE. These problems are solved as follows: first, we establish three kinds of MSE uplink-downlink duality by transforming only the power allocation matrices from uplink channel to downlink channel and vice versa. Second, in the uplink channel, we formulate the power allocation part of each problem ensuring global optimality. Finally, based on the solution of the uplink power allocation and the MSE duality results, for each problem, we propose an iterative algorithm that performs optimization alternatively between the uplink and downlink channels. Computer simulations verify the robustness of the proposed design compared to the non-robust/naive design.

## I. INTRODUCTION

In a multiuser network the uplink-downlink duality approach for solving the downlink optimization problems has received a lot of attention. The achievable sum rate of the broadcast channel (BC) obtained by dirty paper precoding technique has been characterized for multiple-input single-output (MISO) systems [1]. The latter work has been extended in [2] for multiple-input multiple-output (MIMO) systems. These papers analyze the sum rate region of the BC channel by exploiting the duality between BC and multiple access channels (MAC). In [3], the dirty paper rate region has shown to be the capacity region of the Gaussian MIMO BC channel. In [4] and [5], mean-square-error (MSE) based uplink-downlink duality have been exploited. The latter two papers utilize their duality results to solve MSE-based design problems. All of the aforementioned duality are established by assuming that perfect channel state information (CSI) is available at the base station (BS) and mobile stations (MSs). However, due to the inevitability of channel estimation error, CSI can never be perfect. This motivates [6] to establish the MSE duality under imperfect CSI for MISO systems. The latter work is extended in [7] for MIMO case. None of [6] and [7] incorporates antenna correlation in their channel model and neither of these duality can be applied to symbol wise MSE-based problems for MIMO systems. For instance, the duality of [6] and [7] can not be used for the robust symbol wise weighted sum MSE problem. Moreover, while solving the robust sum MSE minimization problem, the authors of [6] and [7] compute $K$ (total number of MSs) scaling factors (see (16) in [6] and [7]) to transfer the total sum average mean-square-error (AMSE) from uplink to downlink channel which is not computationally efficient. As will be seen later in Section IV, we compute only one scaling factor to transfer the sum AMSE from uplink to downlink channel and vice versa. In [8], the MSE uplink-downlink duality has been established by considering imperfect CSI both at the BS and MSs, and with antenna correlation only at the BS. The duality is examined by analyzing the Karush-Kuhn-Tucker conditions for the uplink and downlink channel problems. The latter duality is limited to sum MSE minimization problem.

In [9], we have established three kinds of MSE duality by considering that imperfect CSI is available both at the BS and MSs, and with antenna correlation only at the BS. These duality are established by extending the three level MSE duality of [5] to imperfect CSI. Thus, from the MSE duality perspective, the duality of [9] is more general than that of [6], [7] and [8]. In order to solve general MSE-based robust design problems (see for example (14) in **Case 2** of [9]), the approach of [4] and [10] has been employed where the precoder of each MS is decomposed into a product of unity norm filter and diagonal power allocation matrices, and the decoder of each MS is decomposed into a product of unity norm filter, diagonal scaling factor and the inverse of power allocation matrices (see (15) of [9]). Upon doing so, in [9], we have shown that any MSE-based robust design problem can be solved using alternating optimization framework. From (22) of [9], we have also realized that by employing the same filters and scaling factors in both the uplink and downlink channels, three kinds of AMSE uplink-downlink duality can be established just by transforming the power allocation matrices from uplink channel to downlink channel and vice versa. This motivates us to use the system model shown in Fig. 1. Note that although this system model is known from [4] and [10], the authors of these two papers employ another approach to establish the MSE uplink-downlink duality which is computationally costly.

In the current paper, we consider that the BS and MS antennas exhibit spatial correlations and the CSI at both ends is imperfect. The robustness against imperfect CSI is incorporated into our designs using stochastic approach [8]. In this regard, we first establish three kinds of AMSE duality. Then, as application examples, we consider the joint optimization of transceivers for the following MSE-based robust design problems:

1) The robust minimization of the weighted sum MSE constrained with a total BS power ($\mathcal{P}1$).
2) The robust minimization of the maximum weighted MSE (min-max) constrained with a total BS power ($\mathcal{P}2$).

**Motivations for $\mathcal{P}1$ and $\mathcal{P}2$:** In a multiuser scenario, fairness is an important issue which in general can be achieved by ensuring a minimum level of quality of service (in terms of SINR or MSE) to all users. This applies for both the non-robust and robust designs. The objective function of the robust design problem $\mathcal{P}1$ maintains fairness by allocating the weights in proportion to the priority given to the users. In problem $\mathcal{P}2$, the objective function is the minimization of the maximum weighted MSE which obviously tries to reduce the worst user weighted MSE. Thus, both of the considered problems try to enhance the system performance (with MSE as a performance

The authors would like to thank the Region Wallonne for the financial support of the project MIMOCOM in the framework of which this work has been achieved. Part of this work has been published in the 3rd IEEE International Workshop on CAMSAP, Dec 2009. Tadilo Endeshaw and Luc Vandendorpe are with the Communication and Remote Sensing Laboratory, Université catholique de Louvain, Place du Levant 2, 1348 - Louvain La Neuve, Belgium. Email: {tadilo.bogale, luc.vandendorpe}@uclouvain.be, Phone: +3210478071, Fax: +3210472089. Batu K. Chalise was with the Communication and Remote Sensing Laboratory, Universite catholique de Louvain, Belgium when this work was completed. He is currently with the Center for Advanced Communications, Villanova University, 800 Lancaster Avenue, Villanova, PA, 19085, USA, Phone: +16105197371, Fax: +16105196118, Email: batu.chalise@villanova.edu.

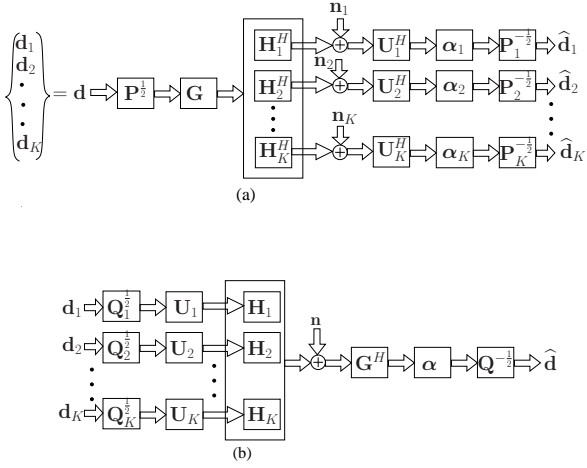

Fig. 1. Multiuser MIMO system. (a) downlink channel. (b) uplink channel.

metric) by taking into account fairness issues.

As $\mathcal{P}1$ and $\mathcal{P}2$ are non-convex, we can not use the convex optimization tools to solve them. Due to this, we first solve the power allocation part of each problem ensuring global optimum. With this solution and the AMSE duality results, like in [9], we propose iterative algorithms for $\mathcal{P}1$ and $\mathcal{P}2$. Thus, this paper has the following contributions:

1) By using the system model shown in Fig. 1, we establish three kinds of AMSE duality known from [5][1] for the aforementioned CSI just by transforming the power allocation matrices from uplink to downlink channel and vice versa. In contrast to the AMSE duality in [5], [6], [7] and [8], our duality can be used to solve all MSE-based problems by using alternating optimization like in [9]. It is worthwhile to mention that one can also extend the duality approach of [4] to imperfect CSI case as the latter duality also requires only the transformation of powers from uplink to downlink channel and vice versa. However, by utilizing our duality, the computational complexity of the latter power transformation can be reduced (this will be clear later in Section V-A). As a consequence, the overall computational cost of alternating optimization algorithm of [4] reduces. Moreover, this work generalizes the hitherto MSE uplink-downlink duality[2].
2) We show that the uplink power allocation part of each problem can be solved ensuring global optimality.
3) Using the uplink power allocation and AMSE duality results, we propose iterative algorithms for $\mathcal{P}1$ and $\mathcal{P}2$.
4) We examine the effects of channel estimation errors and antenna correlations on the system performance.

## II. SYSTEM MODEL

In this section the MIMO downlink and uplink system models are considered. The BS equipped with $N$ transmit antennas is serving $K$ decentralized MSs each having $\{M_k\}_{k=1}^K$ antennas to multiplex $S_k$ symbols. The total number of MS antennas and symbols are $M = \sum_{k=1}^K M_k$ and $S = \sum_{k=1}^K S_k$, respectively. All symbols can be stacked in a data vector $\mathbf{d} = [\mathbf{d}_1^T, \cdots, \mathbf{d}_K^T]^T$, where $\mathbf{d}_k \in C^{S_k \times 1}$ is the symbol vector for the $k$th MS. The MAC channel can be expressed as $\mathbf{H} = [\mathbf{H}_1, \cdots, \mathbf{H}_K]$, where $\mathbf{H}_k^H \in C^{M_k \times N}$ is the channel between the BS and $k$th MS.

Using the system model similar to [10] and as shown in Fig. 1, we collect the transmit powers of all users as $\mathbf{P} = \mathrm{blkdiag}(\mathbf{P}_1, \mathbf{P}_2, \cdots, \mathbf{P}_K)$ and $\mathbf{Q} = \mathrm{blkdiag}(\mathbf{Q}_1, \cdots, \mathbf{Q}_K)$, where $\mathbf{P}_k = \mathrm{diag}(p_{k1}, \cdots, p_{kS_k})$, $\mathbf{Q}_k = \mathrm{diag}(q_{k1}, \cdots, q_{kS_k})$ and $p_{ki}$ ($q_{ki}$) is the downlink (uplink) power allocation for the $i$th symbol of the $k$th user. The overall filter matrix at the BS is $\mathbf{G} = [\mathbf{G}_1, \cdots, \mathbf{G}_K]$, where $\mathbf{G}_k = [\mathbf{g}_{k1} \cdots \mathbf{g}_{kS_k}] \in C^{N \times S_k}$ is the filter matrix for the $k$th user with $\{\mathbf{g}_{ki}^H \mathbf{g}_{ki} = 1\}_{i=1}^{S_k}, k = \{1, \cdots, K\}$. The filters of all users are stacked in a block diagonal matrix $\mathbf{U} = \mathrm{blkdiag}(\mathbf{U}_1, \cdots, \mathbf{U}_K)$, where $\mathbf{U}_k = [\mathbf{u}_{k1} \cdots \mathbf{u}_{kS_k}] \in C^{M_k \times S_k}$ is the filter matrix for the $k$th user with $\{\mathbf{u}_{ki}^H \mathbf{u}_{ki} = 1\}_{i=1}^{S_k}, \forall k$. The scaling factors are accumulated as $\boldsymbol{\alpha} = \mathrm{blkdiag}(\boldsymbol{\alpha}_1, \cdots, \boldsymbol{\alpha}_K)$, where $\boldsymbol{\alpha}_k = \mathrm{diag}(\alpha_{k1}, \cdots, \alpha_{kS_k})$. The entries of $\mathbf{n} = [\mathbf{n}_1^T, \mathbf{n}_2^T, \cdots, \mathbf{n}_K^T]^T$ are assumed to be independent and identically distributed (i.i.d) zero-mean circularly symmetric complex Gaussian (ZMCSCG) random variables all with variance $\sigma^2$. We also assume that $\mathrm{E}\{\mathbf{d}_k \mathbf{d}_k^H\} = \mathbf{I}_{S_k}$, $\mathrm{E}\{\mathbf{d}_k \mathbf{d}_i^H\} = \mathbf{0}, \forall i \neq k$, and $\mathrm{E}\{\mathbf{d}_k \mathbf{n}_k^H\} = \mathbf{0}$.

## III. CHANNEL MODEL

Considering antenna correlation at the BS and MSs, we model the Rayleigh fading MIMO channels as $\mathbf{H}_k^H = \widetilde{\mathbf{R}}_{mk}^{1/2} \mathbf{H}_{wk}^H \mathbf{R}_{bk}^{1/2}, \forall k$, where the elements of $\{\mathbf{H}_{wk}^H\}_{k=1}^K$ are i.i.d ZMCSCG random variables all with unit variance and $\mathbf{R}_{bk} \in C^{N \times N}$, $\widetilde{\mathbf{R}}_{mk} \in C^{M_k \times M_k}$ are antenna correlation matrices at the BS and MSs, respectively [11], [12]. The channel estimation is performed on $\{\mathbf{H}_{wk}^H\}_{k=1}^K$ using an orthogonal training method [13]. Upon doing so, the $k$th user true channel $\mathbf{H}_k^H$ and its minimum-mean-square-error (MMSE) estimate $\widehat{\mathbf{H}}_k^H$ can be related as (see Section II.B of [13] for more details about the channel estimation process)

$$\mathbf{H}_k^H = \widehat{\mathbf{H}}_k^H + \mathbf{R}_{mk}^{1/2} \mathbf{E}_{wk}^H \mathbf{R}_{bk}^{1/2} = \widehat{\mathbf{H}}_k^H + \mathbf{E}_k^H, \ \forall k \quad (1)$$

where $\mathbf{R}_{mk} = (\mathbf{I}_{M_k} + \sigma_{ek}^2 \widetilde{\mathbf{R}}_{mk}^{-1})^{-1}$, $\mathbf{E}_k^H$ is the estimation error and the entries of $\mathbf{E}_{wk}^H$ are i.i.d with $\mathcal{CN}(0, \sigma_{ek}^2)$. In this paper we consider that $\{\mathbf{E}_{wk}^H\}_{k=1}^K$ are unknown but $\{\widehat{\mathbf{H}}_k^H, \mathbf{R}_{bk}, \widetilde{\mathbf{R}}_{mk}$ and $\sigma_{ek}^2\}_{k=1}^K$ are available. We assume that each MS estimates its channel and feeds the estimated channel back to the BS without any error. Thus, both the BS and MSs have the same channel imperfections. The $k$th user estimated signal $\widehat{\mathbf{d}}_k^{DL}$ can be expressed as

$$\widehat{\mathbf{d}}_k^{DL} = \mathbf{P}_k^{-1/2} \boldsymbol{\alpha}_k \mathbf{U}_k^H (\mathbf{H}_k^H \mathbf{G} \mathbf{P}^{1/2} \mathbf{d} + \mathbf{n}_k)$$
$$= \mathbf{P}_k^{-1/2} \boldsymbol{\alpha}_k \mathbf{U}_k^H (\mathbf{H}_k^H \sum_{k=1}^K \mathbf{G}_k \mathbf{P}_k^{1/2} \mathbf{d}_k + \mathbf{n}_k) \quad (2)$$

where $\mathbf{H}_k^H$ is the channel between the BS and the $k$th user, and $\mathbf{n}_k$ is the additive noise at the $k$th MS. The downlink instantaneous MSE matrix of the $k$th user $\boldsymbol{\xi}_k^{DL} = \mathrm{E}_\mathbf{d}\{(\mathbf{d}_k - \widehat{\mathbf{d}}_k^{DL})(\mathbf{d}_k - \widehat{\mathbf{d}}_k^{DL})^H\}$ is given by

$$\boldsymbol{\xi}_k^{DL} = \mathbf{I}_{S_k} + \mathbf{P}_k^{-1/2} \boldsymbol{\alpha}_k \mathbf{U}_k^H (\mathbf{H}_k^H (\sum_{i=1}^K \mathbf{G}_i \mathbf{P}_i \mathbf{G}_i^H) \mathbf{H}_k +$$
$$\sigma^2 \mathbf{I}_{M_k}) \mathbf{U}_k \boldsymbol{\alpha}_k \mathbf{P}_k^{-1/2} - \mathbf{P}_k^{1/2} \mathbf{G}_k^H \mathbf{H}_k \mathbf{U}_k \boldsymbol{\alpha}_k \mathbf{P}_k^{-1/2}$$
$$- \mathbf{P}_k^{-1/2} \boldsymbol{\alpha}_k \mathbf{U}_k^H \mathbf{H}_k^H \mathbf{G}_k \mathbf{P}_k^{1/2}.$$

---

[1]Note: The authors of [5] establish the three kinds of duality by transferring the precoder/decoder pairs from uplink to downlink channel and vice versa.

[2]Note that for the considered CSI model, the MSE uplink-downlink duality can be established using the system model like in [8] and [9]. However, this system model is not convenient to solve general MSE-based robust design problems (for example $\mathcal{P}1$ (**Case 2**) and $\mathcal{P}2$).

Substituting (1) in $\boldsymbol{\xi}_k^{DL}$ and taking the expected value of $\boldsymbol{\xi}_k^{DL}$ over $\mathbf{E}_{w_k}^H$, the downlink AMSEs can be expressed as

$$\overline{\boldsymbol{\xi}}_k^{DL} = \mathrm{E}_{\mathbf{E}_{w_k}^H}\{\boldsymbol{\xi}_k^{DL}\}$$
$$= \mathbf{I} + \mathbf{P}_k^{-1/2}\boldsymbol{\alpha}_k \mathbf{U}_k^H \boldsymbol{\Gamma}_k^{DL}\mathbf{U}_k\boldsymbol{\alpha}_k\mathbf{P}_k^{-1/2} -$$
$$\mathbf{P}_k^{1/2}\mathbf{G}_k^H\widehat{\mathbf{H}}_k\mathbf{U}_k\boldsymbol{\alpha}_k\mathbf{P}_k^{-1/2} - \mathbf{P}_k^{-1/2}\boldsymbol{\alpha}_k\mathbf{U}_k^H\widehat{\mathbf{H}}_k^H\mathbf{G}_k\mathbf{P}_k^{1/2},$$

$$\overline{\xi}_k^{DL} = \mathrm{tr}\{\overline{\boldsymbol{\xi}}_k^{DL}\}$$
$$= S_k + \mathrm{tr}\{\mathbf{P}_k^{-1}\boldsymbol{\alpha}_k^2\mathbf{U}_k^H\boldsymbol{\Gamma}_k^{DL}\mathbf{U}_k - 2\Re\{\mathbf{G}_k^H\widehat{\mathbf{H}}_k\mathbf{U}_k\boldsymbol{\alpha}_k\}\}, \quad (3)$$

$$\overline{\xi}^{DL} = \sum_{k=1}^K \overline{\xi}_k^{DL} \quad (4)$$
$$= S + \sum_{k=1}^K \mathrm{tr}\{\mathbf{P}_k^{-1}\boldsymbol{\alpha}_k^2\mathbf{U}_k^H\boldsymbol{\Gamma}_k^{DL}\mathbf{U}_k - 2\Re\{\mathbf{G}_k^H\widehat{\mathbf{H}}_k\mathbf{U}_k\boldsymbol{\alpha}_k\}\}$$

where $\boldsymbol{\Gamma}_k^{DL} = (\widehat{\mathbf{H}}_k^H\mathbf{GPG}^H\widehat{\mathbf{H}}_k + \sigma_{ek}^2\mathrm{tr}\{\mathbf{R}_{bk}\mathbf{GPG}^H\}\mathbf{R}_{mk} + \sigma^2\mathbf{I}_{M_k})$ and we use the fact $\mathrm{E}_{\mathbf{E}}\{\mathbf{EAE}^H\} = \sigma_e^2\mathrm{tr}\{\mathbf{A}\}\mathbf{I}$, if the entries of $\mathbf{E}$ are i.i.d with $\mathcal{CN}(0, \sigma_e^2)$ and $\mathbf{A}$ is a given matrix. Like in the downlink channel, by defining $\boldsymbol{\Gamma}_c \triangleq [\sum_{i=1}^K(\widehat{\mathbf{H}}_i\mathbf{U}_i\mathbf{Q}_i\mathbf{U}_i^H\widehat{\mathbf{H}}_i^H + \sigma_{ei}^2\mathrm{tr}\{\mathbf{R}_{mi}\mathbf{U}_i\mathbf{Q}_i\mathbf{U}_i^H\}\mathbf{R}_{bi}) + \sigma^2\mathbf{I}_N]$, the uplink channel AMSEs are given by

$$\overline{\boldsymbol{\xi}}_k^{UL} = \mathbf{I}_{S_k} + \mathbf{Q}_k^{-1/2}\boldsymbol{\alpha}_k\mathbf{G}_k^H\boldsymbol{\Gamma}_c\mathbf{G}_k\boldsymbol{\alpha}_k\mathbf{Q}_k^{-1/2} - \mathbf{Q}_k^{-1/2}\boldsymbol{\alpha}_k\cdot$$
$$\mathbf{G}_k^H\widehat{\mathbf{H}}_k\mathbf{U}_k\mathbf{Q}_k^{1/2} - \mathbf{Q}_k^{1/2}\mathbf{U}_k^H\widehat{\mathbf{H}}_k^H\mathbf{G}_k\boldsymbol{\alpha}_k\mathbf{Q}_k^{-1/2}. \quad (5)$$

$$\overline{\xi}_k^{UL} = \mathrm{tr}\{\overline{\boldsymbol{\xi}}_k^{UL}\}$$
$$= S_k + \mathrm{tr}\{\mathbf{Q}_k^{-1}\boldsymbol{\alpha}_k^2\mathbf{G}_k^H\boldsymbol{\Gamma}_c\mathbf{G}_k - 2\Re\{\boldsymbol{\alpha}_k\mathbf{G}_k^H\widehat{\mathbf{H}}_k\mathbf{U}_k\}\}. \quad (6)$$

$$\overline{\xi}^{UL} = \sum_{k=1}^K \overline{\xi}_k^{UL}$$
$$= S + \sum_{k=1}^K \mathrm{tr}\{\mathbf{Q}_k^{-1}\boldsymbol{\alpha}_k^2\mathbf{G}_k^H\boldsymbol{\Gamma}_c\mathbf{G}_k - 2\Re\{\boldsymbol{\alpha}_k\mathbf{G}_k^H\widehat{\mathbf{H}}_k\mathbf{U}_k\}\}. \quad (7)$$

## IV. AVERAGE MEAN SQUARE ERROR UPLINK-DOWNLINK DUALITY

As we mentioned in Section I, our AMSE duality generalizes the work of [9] to the case where the BS and MS antennas are spatially correlated, and both the BS and MSs have imperfect CSI. Thus, in this section, we transfer the sum AMSE, user wise AMSE and symbol wise AMSEs from the uplink to downlink channel and vice versa.

### A. AMSE transfer from uplink to downlink channel

*1) Total sum AMSE transfer:* For a given uplink sum AMSE (with a transmit power $\mathbf{Q}$), we can achieve the same sum AMSE in the downlink channel by using a positive $\beta$ which satisfies $\mathbf{P} = \beta\boldsymbol{\alpha}^2\mathbf{Q}^{-1}$. Substituting $\mathbf{P}$ in (4), equating $\overline{\xi}^{DL} = \overline{\xi}^{UL}$ and after some simple derivations, $\beta$ can be determined as

$$\beta = \mathrm{tr}\{\mathbf{Q}\}/\mathrm{tr}\{\mathbf{Q}^{-1}\boldsymbol{\alpha}^2\}. \quad (8)$$

As can be seen from (8), the scaling factor $\beta$ does not depend on $\{\sigma_{ek}^2\}_{k=1}^K$. This can be seen from (4) and (7), after substituting $\{\boldsymbol{\Gamma}_k^{DL}\}_{k=1}^K$ and $\boldsymbol{\Gamma}_c$, where $\{\sigma_{ek}^2\}_{k=1}^K$ are amplified by the same factor. The downlink power is given by $P_{sum}^{DL} = \mathrm{tr}\{\mathbf{P}\} = \mathrm{tr}\{\beta\boldsymbol{\alpha}^2\mathbf{Q}^{-1}\} = \mathrm{tr}\{\mathbf{Q}\} = P_{sum}^{UL}$. Thus, the same sum power is allocated in both channels.

*2) User wise AMSE transfer:* Given the $k$th user AMSE in the uplink channel with $\{\mathbf{Q}_k\}_{k=1}^K \neq \mathbf{0}$, this user can achieve the same AMSE in the downlink channel if $\mathbf{P}_k$ is computed by

$$\mathbf{P}_k = \beta_k\boldsymbol{\alpha}_k^2\mathbf{Q}_k^{-1}. \quad (9)$$

Substituting (9) in (3), then equating $\overline{\xi}_k^{UL} = \overline{\xi}_k^{DL}$ and after some mathematical manipulations we obtain

$$\beta_k\sigma^2\mathrm{tr}\{\mathbf{Q}_k^{-1}\boldsymbol{\alpha}_k^2\} + \sum_{i=1,i\neq k}^K \beta_k(\|\mathbf{Q}_k^{-1/2}\boldsymbol{\alpha}_k\mathbf{G}_k^H\widehat{\mathbf{H}}_i\mathbf{U}_i\mathbf{Q}_i^{1/2}\|_F^2 +$$
$$\sigma_{ei}^2\mathrm{tr}\{\mathbf{R}_{mi}\mathbf{U}_i\mathbf{Q}_i\mathbf{U}_i^H\}\|\mathbf{R}_{bi}^{1/2}\mathbf{G}_k\boldsymbol{\alpha}_k\mathbf{Q}_k^{-1/2}\|_F^2) =$$
$$\sigma^2\mathrm{tr}\{\mathbf{Q}_k\} + \sum_{i=1,i\neq k}^K \beta_i(\|\boldsymbol{\alpha}_i\mathbf{Q}_i^{-1/2}\mathbf{G}_i^H\widehat{\mathbf{H}}_k\mathbf{U}_k\mathbf{Q}_k^{1/2}\|_F^2 +$$
$$\sigma_{ek}^2\mathrm{tr}\{\mathbf{R}_{mk}\mathbf{U}_k\mathbf{Q}_k\mathbf{U}_k^H\}\|\mathbf{R}_{bk}^{1/2}\mathbf{G}_i\boldsymbol{\alpha}_i\mathbf{Q}_i^{-1/2}\|_F^2). \quad (10)$$

After applying (10) for all users, we can form the following system of linear equations

$$\mathbf{X} \cdot [\beta_1, \ldots, \beta_K]^T = \sigma^2 [\mathrm{tr}\{\mathbf{Q}_1\}, \ldots, \mathrm{tr}\{\mathbf{Q}_K\}]^T \quad (11)$$

where $[\mathbf{X}]_{k,l} = \quad (12)$

$$\begin{cases} \sigma^2\mathrm{tr}\{\mathbf{Q}_k^{-1}\boldsymbol{\alpha}_k^2\} + \sum_{i=1,i\neq k}^K(\|\mathbf{Q}_k^{-1/2}\boldsymbol{\alpha}_k\mathbf{G}_k^H\widehat{\mathbf{H}}_i\mathbf{U}_i\mathbf{Q}_i^{1/2}\|_F^2 + \\ \sigma_{ei}^2\|\mathbf{R}_{mi}^{1/2}\mathbf{U}_i\mathbf{Q}_i^{1/2}\|_F^2\|\mathbf{R}_{bi}^{1/2}\mathbf{G}_k\boldsymbol{\alpha}_k\mathbf{Q}_k^{-1/2}\|_F^2), & k = l \\ -(\|\boldsymbol{\alpha}_l\mathbf{Q}_l^{-1/2}\mathbf{G}_l^H\widehat{\mathbf{H}}_k\mathbf{U}_k\mathbf{Q}_k^{1/2}\|_F^2 + \\ \sigma_{ek}^2\|\mathbf{R}_{mk}^{1/2}\mathbf{U}_k\mathbf{Q}_k^{1/2}\|_F^2\|\mathbf{R}_{bk}^{1/2}\mathbf{G}_l\boldsymbol{\alpha}_l\mathbf{Q}_l^{-1/2}\|_F^2), & k \neq l. \end{cases}$$

It can be shown that if $\sigma^2 > 0$ then $\{\beta_k\}_{k=1}^K$ of (11) are strictly positive [5], [9]. Thus, the $k$th user AMSE can be transferred from uplink to downlink channel. Summing up the left-hand and right-hand sides of (11) and cancelling $\sigma^2$ in both sides yields $P_{sum}^{DL} = \sum_{i=1}^K \beta_i \mathrm{tr}\{\mathbf{Q}_i^{-1}\boldsymbol{\alpha}_i^2\} = \sum_{i=1}^K \mathrm{tr}\{\mathbf{P}_i\} = \sum_{i=1}^K \mathrm{tr}\{\mathbf{Q}_i\} = P_{sum}^{UL}$. Thus, the same sum power is allocated in both the uplink and downlink channels.

### B. AMSE transfer from downlink to uplink channel

To complete the duality, in this section we examine the AMSE transfer from the downlink to uplink channel.

*1) Total sum AMSE transfer:* Similar to Subsection IV-A.1, the sum AMSE can be transferred from the downlink to uplink channel by using a nonzero scaling factor $\widetilde{\beta}$ which satisfies $\mathbf{Q} = \widetilde{\beta}\boldsymbol{\alpha}^2\mathbf{P}^{-1}$. Substituting $\mathbf{Q}$ in (7) and then equating $\overline{\xi}^{UL} = \overline{\xi}^{DL}$, $\widetilde{\beta}$ is determined as

$$\widetilde{\beta} = \mathrm{tr}\{\mathbf{P}\}/\mathrm{tr}\{\mathbf{P}^{-1}\boldsymbol{\alpha}^2\}. \quad (13)$$

*2) User wise AMSE transfer:* Given the $k$th user downlink AMSE with $\{\mathbf{P}_k\}_{k=1}^K \neq \mathbf{0}$, this user can achieve the same AMSE in the uplink channel if $\mathbf{Q}_k$ is computed by $\mathbf{Q}_k = \widetilde{\beta}_k\boldsymbol{\alpha}_k^2\mathbf{P}_k^{-1}$. Like in Subsection IV-A.2, by substituting $\mathbf{Q}_k$ in (6), equating $\overline{\xi}_k^{UL} = \overline{\xi}_k^{DL}$ and after some mathematical manipulations, the scaling factors $\{\widetilde{\beta}_k\}_{k=1}^K$ are determined by solving the following system of linear equations.

$$\mathbf{T} \cdot [\widetilde{\beta}_1, \ldots, \widetilde{\beta}_K]^T = \sigma^2 [\mathrm{tr}\{\mathbf{P}_1\}, \ldots, \mathrm{tr}\{\mathbf{P}_K\}]^T \quad (14)$$

where $[\mathbf{T}]_{k,l} = \quad (15)$

$$\begin{cases} \sigma^2\mathrm{tr}\{\mathbf{P}_k^{-1}\boldsymbol{\alpha}_k^2\} + \sum_{i=1,i\neq k}^K(\|\mathbf{P}_k^{-1/2}\boldsymbol{\alpha}_k\mathbf{U}_k^H\widehat{\mathbf{H}}_k^H\mathbf{G}_i\mathbf{P}_i^{1/2}\|_F^2 + \\ \sigma_{ek}^2\|\mathbf{R}_{bk}^{1/2}\mathbf{G}_i\mathbf{P}_i^{1/2}\|_F^2\|\mathbf{R}_{mk}^{1/2}\mathbf{U}_k\mathbf{P}_k^{-1/2}\boldsymbol{\alpha}_k\|_F^2), & k = l \\ -(\|\mathbf{P}_l^{-1/2}\boldsymbol{\alpha}_l\mathbf{U}_l^H\widehat{\mathbf{H}}_l^H\mathbf{G}_k\mathbf{P}_k^{1/2}\|_F^2 + \\ \sigma_{el}^2\|\mathbf{R}_{ml}^{1/2}\mathbf{U}_l\boldsymbol{\alpha}_l\mathbf{P}_l^{-1/2}\|_F^2\|\mathbf{R}_{bl}^{1/2}\mathbf{G}_k\mathbf{P}_k^{1/2}\|_F^2), & k \neq l. \end{cases}$$

The symbol wise AMSE transfer (from the uplink channel to downlink channel and vice versa) can be examined similar to Subsections IV-A.2 and IV-B.2. The details are omitted due to space constraint.

## V. APPLICATIONS OF AMSE DUALITY

To show the applications of our AMSE duality, in this section, we examine the problem of jointly designing the precoders and decoders for the downlink multiuser MIMO systems to minimize: (i) the weighted sum MSE under a total BS power constraint ($\mathcal{P}1$) and (ii) the maximum weighted user AMSE constrained with a total BS power ($\mathcal{P}2$). Both design problems provide robustness against the channel uncertainties.

### A. The robust weighted sum MSE minimization problem ($\mathcal{P}1$)

In the downlink channel, the robust weighted sum MSE minimization problem ($\mathcal{P}1$) can be formulated by

$$\min_{\mathbf{G}_k, \mathbf{U}_k, \boldsymbol{\alpha}_k, \mathbf{P}_k} \mathrm{E}_{\mathbf{E}_{wk}^H} \sum_{k=1}^{K} \tau_k \mathrm{tr}\{\boldsymbol{\xi}_k^{DL}\} = \sum_{k=1}^{K} \tau_k \overline{\xi}_k^{DL}$$

$$\mathrm{s.t} \quad \sum_{k=1}^{K} \mathrm{tr}\{\mathbf{P}_k\} \leq P_{max},$$

$$\{\mathbf{g}_{ki}^H \mathbf{g}_{ki} = \mathbf{u}_{ki}^H \mathbf{u}_{ki} = 1, p_{ki} > 0\}_{i=1}^{S_k}, \forall k \quad (16)$$

where $\tau_k$ is the AMSE weighting factor of the $k$th user (when $\{\tau_k = 1\}_{k=1}^{K}$, (16) simplifies to sum AMSE problem). In $\sum_{k=1}^{K} \tau_k \overline{\xi}_k^{DL}$, since the precoders of all users are coupled, $\mathcal{P}1$ has more complicated mathematical structure than its dual uplink problem [4], [5]. Due to this, we examine the dual uplink problem of (16) which is expressed as

$$\min_{\{\mathbf{G}_k, \mathbf{U}_k, \boldsymbol{\alpha}_k, \mathbf{Q}_k\}_{k=1}^{K}} \sum_{k=1}^{K} \tau_k \overline{\xi}_k^{UL}$$

$$\mathrm{s.t} \quad \sum_{k=1}^{K} \mathrm{tr}\{\mathbf{Q}_k\} \leq P_{\max}$$

$$\{\mathbf{g}_{ki}^H \mathbf{g}_{ki} = \mathbf{u}_{ki}^H \mathbf{u}_{ki} = 1, q_{ki} > 0\}_{i=1}^{S_k}, \forall k. \quad (17)$$

For convenience, we consider (17) for the following two cases.

**Case 1:** When $\{\tau_k = 1, \widehat{\mathbf{R}}_{mk} = \mathbf{I}_{M_k}, \mathbf{R}_{bk} = \mathbf{R}_b \text{ and } \sigma_{ek}^2 = \sigma_e^2\}_{k=1}^{K}$: In this case, first, for a fixed transmitter, the $k$th user receiver is optimized by using the minimum average mean-square-error (MAMSE) method which yields

$$\widetilde{\mathbf{G}}_k \triangleq \mathbf{G}_k \boldsymbol{\alpha}_k \mathbf{Q}_k^{-1/2} = \mathbf{\Gamma} \widehat{\mathbf{H}}_k \mathbf{U}_k \mathbf{Q}_k^{1/2} \quad (18)$$

where $\mathbf{\Gamma} = [\sum_{i=1}^{K}(\widehat{\mathbf{H}}_i \mathbf{U}_i \mathbf{Q}_i \mathbf{U}_i^H \widehat{\mathbf{H}}_i^H + \widetilde{\sigma}_e^2 \mathrm{tr}\{\mathbf{Q}_i\} \mathbf{R}_b) + \sigma^2 \mathbf{I}_N]^{-1}$ and $\widetilde{\sigma}_e^2 = \sigma_e^2/(\sigma_e^2+1)$. Then, after substituting $\widetilde{\mathbf{G}}_k$ in $\overline{\xi}_k^{UL}$, we get the $k$th user MAMSE matrix as $\widetilde{\boldsymbol{\xi}}_k^{UL} = \mathbf{I}_{S_k} - \mathbf{Q}_k^{1/2} \mathbf{U}_k^H \widehat{\mathbf{H}}_k^H \mathbf{\Gamma} \widehat{\mathbf{H}}_k \mathbf{U}_k \mathbf{Q}_k^{1/2}$. It follows that

$$\sum_{k=1}^{K} \mathrm{tr}\{\widetilde{\boldsymbol{\xi}}_k^{UL}\} = \sum_{k=1}^{K} \mathrm{tr}\{\mathbf{I}_{S_k} - \mathbf{Q}_k^{1/2} \mathbf{U}_k^H \widehat{\mathbf{H}}_k^H \mathbf{\Gamma} \widehat{\mathbf{H}}_k \mathbf{U}_k \mathbf{Q}_k^{1/2}\} \quad (19)$$

$$= S - N + \mathrm{tr}\{(\sum_{k=1}^{K} \widetilde{\sigma}_e^2 \mathrm{tr}\{\mathbf{Q}_k\} \mathbf{R}_b + \sigma^2 \mathbf{I}_N) \mathbf{\Gamma}\}$$

where the second equality is derived using the matrix inversion Lemma and the fact $(\mathbf{AB})^{-1} = \mathbf{B}^{-1} \mathbf{A}^{-1}$ [14]. Thus, (17) can be solved by applying a two step approach. First $\mathbf{U}_k$ and $\mathbf{Q}_k$ are optimized by

$$\min_{\{\mathbf{U}_k, \mathbf{Q}_k\}_{k=1}^{K}} \mathrm{tr}\{(\sum_{k=1}^{K} \widetilde{\sigma}_e^2 \mathrm{tr}\{\mathbf{Q}_k\} \mathbf{R}_b + \sigma^2 \mathbf{I}_N) \mathbf{\Gamma}\}$$

$$\mathrm{s.t} \quad \sum_{k=1}^{K} \mathrm{tr}\{\mathbf{Q}_k\} \leq P_{max}, \{\mathbf{u}_{ki}^H \mathbf{u}_{ki} = 1, q_{ki} > 0\}_{i=1}^{S_k}, \forall k, \quad (20)$$

and then the optimum $\mathbf{G}_k$ and $\boldsymbol{\alpha}_k$ are computed by using (18). Note that in (18), $\mathbf{G}_k$ and $\boldsymbol{\alpha}_k$ are obtained such that $\{\mathbf{g}_{ki}^H \mathbf{g}_{ki} = 1\}_{i=1}^{S_k}$, $\forall k$ and $\{\boldsymbol{\alpha}_k\}_{k=1}^{K}$ are diagonal matrices.

Using matrix inversion Lemma, (19) can also be written in terms of $\overline{\mathbf{Q}} \triangleq \mathbf{Q}/\mathrm{tr}\{\mathbf{Q}\}$ as $\sum_{k=1}^{K} \mathrm{tr}\{\widetilde{\boldsymbol{\xi}}_k^{UL}\} = \mathrm{tr}\{(\mathbf{I}_S + \overline{\mathbf{Q}}^{1/2} \mathbf{U}^H \widehat{\mathbf{H}}^H (\widetilde{\sigma}_e^2 \mathbf{R}_b + \frac{\sigma^2}{\mathrm{tr}\{\mathbf{Q}\}} \mathbf{I}_N)^{-1} \widehat{\mathbf{H}} \mathbf{U} \overline{\mathbf{Q}}^{1/2})^{-1}\}$. According to [15], for the given $\overline{\mathbf{Q}}$, $\sum_{k=1}^{K} \mathrm{tr}\{\widetilde{\boldsymbol{\xi}}_k^{UL}\}$ is a non-increasing function of $\mathrm{tr}\{\mathbf{Q}\} = P_{\mathrm{sum}}$. Since the difference between (19) and the objective function of (20) is only the constant term $S - N$, it is clear that the latter objective function is also non-increasing in $P_{\mathrm{sum}}$. By defining $\{\overline{\mathbf{U}}_k \triangleq \mathbf{U}_k \mathbf{Q}_k \mathbf{U}_k^H\}_{k=1}^{K}$, problem (20) can thus be equivalently formulated as

$$\min_{\{\overline{\mathbf{U}}_k\}_{k=1}^{K}} \mathrm{tr}\{(\widetilde{\sigma}_e^2 P_{\max} \mathbf{R}_b + \sigma^2 \mathbf{I}_N) \widetilde{\mathbf{\Gamma}}\}$$

$$\mathrm{s.t} \quad \sum_{i=1}^{K} \mathrm{tr}\{\overline{\mathbf{U}}_k\} = P_{max}, \overline{\mathbf{U}}_k \succeq 0,$$

$$\mathrm{rank}\{\overline{\mathbf{U}}_k\} = \min(M_k, S_k), \forall k \quad (21)$$

where $\widetilde{\mathbf{\Gamma}} = [\sum_{i=1}^{K} \widehat{\mathbf{H}}_i \overline{\mathbf{U}}_i \widehat{\mathbf{H}}_i^H + \widetilde{\sigma}_e^2 P_{max} \mathbf{R}_b + \sigma^2 \mathbf{I}_N]^{-1}$. If we ignore (relax) the rank-constraint of (21), the above problem can be formulated as a semi-definite programming (SDP) problem for which global optimum is guaranteed [16], [17], [18]. Now, if the optimal solution of this SDP satisfies $\mathrm{rank}\{\overline{\mathbf{U}}_k\} = \min(M_k, S_k)$, the latter solution can be considered as a global minimizer of (21), otherwise, the solution is deemed as the lower bound solution of (21). After computing the solution of (21), the optimum $\{\mathbf{U}_k, \mathbf{Q}_k\}_{k=1}^{K}$ are determined from the eigenvalue decomposition of $\{\overline{\mathbf{U}}_k\}_{k=1}^{K}$ (see Table I of [4]). It turns out that the optimum (either local or global) solution of (16) is computed by using our sum AMSE transfer (see Section IV-A.1).

For $\{M_k = S_k = L\}_{k=1}^{K}$, the approach of [4] requires $O(K^3 N^3)$ operations to transfer the powers from uplink to downlink channel (see appendix of [4]) whereas our proposed method needs only $O(KL)$ operations. Thus, as claimed in Section I, the proposed power transformation requires less computation than that of in [4].

**Case 2:** For any $\{\tau_k, \widehat{\mathbf{R}}_{mk}, \mathbf{R}_{bk} \text{ and } \sigma_{ek}^2\}_{k=1}^{K}$: In such general case, (17) can not be formulated as an SDP problem. Thus, the solution method discussed for **Case 1** can not be applied. Due to this, here we first formulate the power allocation part of (17) as a Geometric Programming (GP) for which global optimality is guaranteed. Then, based on the solution of GP, MAMSE receiver and AMSE duality results, we solve (16) using the alternating optimization method like in [9]. To this end, we rewrite $\overline{\xi}_k^{UL}$ into a form which is suitable for the GP formulation. Using (6), we can express $\overline{\xi}_k^{UL}$ as

$$\overline{\xi}_k^{UL} = \lambda_k + \widetilde{q}_k^{-1} \sum_{i=1, i \neq k}^{K} \widetilde{q}_i \upsilon_{ki} + \sigma^2 \widetilde{q}_k^{-1} \vartheta_k \quad (22)$$

where $\mathbf{Q}_k = \widetilde{q}_k \widetilde{\mathbf{Q}}_k$, $\mathrm{tr}\{\widetilde{\mathbf{Q}}_k\} = 1$, $\widetilde{\mathbf{U}}_k = \mathbf{U}_k \widetilde{\mathbf{Q}}_k \mathbf{U}_k^H$, $\lambda_k = \mathrm{tr}\{\widetilde{\mathbf{Q}}_k^{-1} \boldsymbol{\alpha}_k^2 \mathbf{G}_k^H (\widehat{\mathbf{H}}_k \widetilde{\mathbf{U}}_k \widehat{\mathbf{H}}_k^H + \sigma_{ek}^2 \mathrm{tr}\{\mathbf{R}_{mk} \widetilde{\mathbf{U}}_k\} \mathbf{R}_{bk}) \mathbf{G}_k - \boldsymbol{\alpha}_k \mathbf{G}_k^H \widehat{\mathbf{H}}_k \mathbf{U}_k - \mathbf{U}_k^H \widehat{\mathbf{H}}_k^H \mathbf{G}_k \boldsymbol{\alpha}_k\} + S_k$, $\vartheta_k = \mathrm{tr}\{\widetilde{\mathbf{Q}}_k^{-1} \boldsymbol{\alpha}_k^2\}$ and $\upsilon_{ki} = \mathrm{tr}\{\widetilde{\mathbf{Q}}_k^{-1} \boldsymbol{\alpha}_k^2 \mathbf{G}_k^H (\widehat{\mathbf{H}}_i \widetilde{\mathbf{U}}_i \widehat{\mathbf{H}}_i^H + \sigma_{ei}^2 \mathrm{tr}\{\mathbf{R}_{mi} \widetilde{\mathbf{U}}_i\} \mathbf{R}_{bi}) \mathbf{G}_k\}$. Once again, we can rewrite $\overline{\xi}_k^{UL}$ as

$$\overline{\xi}_k^{UL} = \widetilde{q}_k^{-1} [\mathbf{Y} \widetilde{\mathbf{q}} + \sigma^2 \boldsymbol{\vartheta}]_k \quad (23)$$

where $\widetilde{\mathbf{q}} = [\widetilde{q}_1, \cdots, \widetilde{q}_K]$, $\boldsymbol{\vartheta} = [\vartheta_1, \cdots, \vartheta_K]^T$ and

$$[\mathbf{Y}]_{k,i} = \begin{cases} \lambda_k, & \text{for } k = i \\ \upsilon_{k,i}, & \text{for } k \neq i. \end{cases}$$

As can be seen from the above equation, for fixed $\{\widetilde{\mathbf{Q}}_k, \widetilde{\mathbf{U}}_k$ and $\boldsymbol{\alpha}_k\}_{k=1}^K$, (23) is a posynomial (where $\widetilde{\mathbf{q}} = [\widetilde{q}_1, \cdots, \widetilde{q}_K]$ are the variables). Thus, the power allocation part of (17) is formulated by the following GP problem.

$$\min_{\{\widetilde{q}_k\}_{k=1}^K} \sum_{k=1}^K \tau_k \overline{\xi}_k^{UL}, \text{ s.t } \|\widetilde{\mathbf{q}}\|_1 \leq P_{max}, \widetilde{q}_k > 0, \forall k. \quad (24)$$

It is important to observe that when perfect CSI is available at the transmitter and receivers, the power allocation parts of rate-based optimization problems have been formulated as GPs [10]. Moreover, for the perfect CSI case, the authors of [10] have shown the connection between rate-based and MSE-based optimization problems. However, to the best of our knowledge, the relation between robust rate-based and robust MSE-based optimization problems is not known (under the stochastic robust design approach). We believe that if the CSI is imperfect both at the BS and MSs, the derivation of rate expression is much involved. Even for the case with perfect CSI at the receiver and imperfect CSI at the transmitter, the stochastic sum rate robust design problem is not easy to solve exactly, and requires a number of bounds and approximations. Such a robust design problem involves expectation of a logarithmic term containing an inverse matrix. Thus, the extension of [10] to robust design case is still an open problem.

Using the solution of (24) and the user wise AMSE duality results, we solve (16) by using the alternating optimization technique similar to that of [9]. In general, we can optimize the powers and filters in many possible orders. In this paper we present a particular algorithm where optimization is started in the uplink channel.

*1) Uplink channel:* In the uplink channel first (24) is solved. With the solution $\{\widetilde{q}_k\}_{k=1}^K$, the powers $\{\mathbf{Q}_k = \widetilde{q}_k \widetilde{\mathbf{Q}}_k\}_{k=1}^K$ are computed and then $\{\mathbf{G}_k$ and $\boldsymbol{\alpha}_k\}_{k=1}^K$ are updated by the following uplink MAMSE receiver

$$\mathbf{G}_k \boldsymbol{\alpha}_k = \boldsymbol{\Gamma}_c^{-1} \widehat{\mathbf{H}}_k \mathbf{U}_k \mathbf{Q}_k. \quad (25)$$

*2) Downlink channel:* Now we switch the optimization to the downlink channel. Thus, we first ensure the same performance as in the uplink channel ($\{\overline{\xi}_k^{DL_1} = \overline{\xi}_k^{UL}\}_{k=1}^K$) by using our user wise AMSE transfer (9). Then, for fixed $\{\mathbf{P}_k\}_{k=1}^K$, the matrices $\{\mathbf{U}_k$ and $\boldsymbol{\alpha}_k\}_{k=1}^K$ are updated by the downlink MAMSE receiver which is given as

$$\mathbf{U}_k \boldsymbol{\alpha}_k = (\boldsymbol{\Gamma}_k^{DL})^{-1} \widehat{\mathbf{H}}_k^H \mathbf{G}_k \mathbf{P}_k. \quad (26)$$

At this stage, the $k$th user achieves a new AMSE $\triangleq \overline{\xi}_k^{DL_2} \leq \overline{\xi}_k^{DL_1}$.

*3) Uplink channel:* Like in Step (2) above, we first ensure the same performance as in the downlink channel ($\{\overline{\xi}_k^{UL_1} = \overline{\xi}_k^{DL_2}\}_{k=1}^K$) and then we update $\{\mathbf{G}_k$ and $\boldsymbol{\alpha}_k\}_{k=1}^K$ by (25). We observe less overall computational time if the latter two steps are performed before proceeding to the next iteration. The detailed iterative steps to solve (16) are summarized in **Table I** (**Algorithm I**).

### B. The robust weighted MSE min-max problem ($\mathcal{P}2$)

In the downlink channel, for given user wise AMSE weights $\{\eta_k\}_{k=1}^K$, $\mathcal{P}2$ can be formulated by

$$\min_{\mathbf{P}_k, \mathbf{G}_k, \mathbf{U}_k, \boldsymbol{\alpha}_k} \max \frac{\mathrm{E}_{\mathbf{E}_{wk}^H} \mathrm{tr}\{\xi_k^{DL}\}}{\eta_k} = \frac{\{\overline{\xi}_k^{DL}\}}{\eta_k}$$

$$\text{s.t } \sum_{k=1}^K \mathrm{tr}\{\mathbf{P}_k\} \leq P_{max}$$

$$\{\mathbf{g}_{ki}^H \mathbf{g}_{ki} = \mathbf{u}_{ki}^H \mathbf{u}_{ki} = 1, p_{ki} > 0\}_{i=1}^{S_k}, \forall k. \quad (27)$$

Here we first solve the power allocation part of (27), then we use the solution framework of $\mathcal{P}1$ (**Case 2**) to jointly optimize the transceivers. To this end, for fixed $\{\widetilde{\mathbf{Q}}_k, \boldsymbol{\alpha}_k, \mathbf{G}_k\}_{k=1}^K$, the uplink power allocation part of (27) reads

$$\mu^{UL} \triangleq \min_{\widetilde{q}_k} \max \frac{\overline{\xi}_k^{UL}}{\eta_k}, \text{ s.t } \|\widetilde{\mathbf{q}}\|_1 \leq P_{max}, \widetilde{q}_k > 0, \forall k. \quad (28)$$

The global optimal solution of the above optimization problem satisfies the following relations [19]

$$\mu^{UL} = \frac{\overline{\xi}_k^{UL}}{\eta_k}, \forall k \text{ and } \|\widetilde{\mathbf{q}}\|_1 = P_{max}. \quad (29)$$

Moreover, by defining $\boldsymbol{\eta} \triangleq \mathrm{diag}\{\eta_1, \eta_2, \cdots, \eta_K\}$ the following eigensystem can be formed from (23) and (29).

$$\boldsymbol{\Omega} \begin{bmatrix} \widetilde{\mathbf{q}} \\ 1 \end{bmatrix} = \mu^{UL} \begin{bmatrix} \widetilde{\mathbf{q}} \\ 1 \end{bmatrix} \quad (30)$$

where

$$\boldsymbol{\Omega} = \begin{bmatrix} \boldsymbol{\eta}^{-1} \mathbf{Y} & \sigma^2 \boldsymbol{\eta}^{-1} \boldsymbol{\vartheta} \\ \frac{1}{P_{max}} \mathbf{1}_K^T \boldsymbol{\eta}^{-1} \mathbf{Y} & \frac{\sigma^2}{P_{max}} \mathbf{1}_K^T \boldsymbol{\eta}^{-1} \boldsymbol{\vartheta} \end{bmatrix}. \quad (31)$$

Therefore, the optimal solution of (28) is given by $\mu^{UL} = \lambda_{max}(\boldsymbol{\Omega})$ and $[\widetilde{\mathbf{q}} \; 1]^T$ is the eigenvector of $\boldsymbol{\Omega}$ corresponding to $\mu^{UL}$ [19]. By using the optimal $\widetilde{\mathbf{q}}$ of (28), MAMSE receiver and AMSE duality results, (27) can be solved as shown in **Table I** (**Algorithm II**).

**TABLE I**: Iterative solution for problems $\mathcal{P}1$ (16) and $\mathcal{P}2$ (27)
Initialization: Set equal powers for all symbols, i.e., $\{\mathbf{Q}_k = \frac{P_{max}}{S} \mathbf{I}_{S_k}\}_{k=1}^K$ and $\mathbf{U}_k \in \mathcal{C}^{M_k \times S_k}$ as the first $S_k$ right-hand singular value decomposition vectors of $\widehat{\mathbf{H}}_k$, $\forall k$ and then update $\{\mathbf{G}_k$ and $\boldsymbol{\alpha}_k\}_{k=1}^K$ by (25).

    repeat    **Virtual uplink channel**
1) Set $\{\widetilde{\mathbf{Q}}_k = \mathbf{Q}_k/\mathrm{tr}\{\mathbf{Q}_k\}\}_{k=1}^K$.
2) For $\mathcal{P}1$, compute $\{\widetilde{q}_k\}_{k=1}^K$ using (24) (**Algorithm I**).
3) For $\mathcal{P}2$, with the given $\boldsymbol{\eta} = \mathrm{diag}\{\eta_1, \eta_2, \cdots, \eta_K\}$, compute $\{\widetilde{q}_k\}_{k=1}^K$ and $\mu^{UL}$ using (31). In this step the power constraint $\|\widetilde{\mathbf{q}}\|_1 = P_{max}$ is ensured by scaling the eigenvector corresponding to $\lambda_{max}(\boldsymbol{\Omega})$ such that the last element equals 1. (**Algorithm II**)
4) Update $\{\mathbf{Q}_k = \widetilde{q}_k \widetilde{\mathbf{Q}}_k\}_{k=1}^K$. Using the latter $\{\mathbf{Q}_k\}_{k=1}^K$, update $\{\mathbf{G}_k$ and $\boldsymbol{\alpha}_k\}_{k=1}^K$ by (25). Then, compute $\{\beta_k\}_{k=1}^K$ with (11).
    **Downlink Channel**
5) Set $\{\mathbf{P}_k = \beta_k \boldsymbol{\alpha}_k^2 \mathbf{Q}_k^{-1}\}_{k=1}^K$. Using this $\{\mathbf{P}_k\}_{k=1}^K$, update $\{\mathbf{U}_k$ and $\boldsymbol{\alpha}_k\}_{k=1}^K$ by (26). Then, compute $\{\widetilde{\beta}_k\}_{k=1}^K$ with (14).
    **Virtual uplink Channel**.
6) Set $\{\mathbf{Q}_k = \widetilde{\beta}_k \boldsymbol{\alpha}_k^2 \mathbf{P}_k^{-1}\}$. Using the latter $\{\mathbf{Q}_k\}_{k=1}^K$, update $\{\mathbf{G}_k$ and $\boldsymbol{\alpha}_k\}_{k=1}^K$ by (25).
    until convergence

**Convergence**: It can be shown that **Algorithms I** and **II** are convergent [10], [19]. Different initializations affect the convergence speed of both algorithms. In most of our simulations (> 95%), we observe fast convergence when the initialization is performed as in this table.
**Global optimality**: Since $\mathcal{P}1$ (**Case 2**) and $\mathcal{P}2$ are non-convex, we can not prove the global optimality of **Algorithms I** and **II**. However, for $\mathcal{P}1$ (**Case 1** with $M_k = S_k$), simulation results show that **Algorithm I** achieves global minimum (see the next Section).

## VI. SIMULATION RESULTS

In all simulations, we take $K = 2$, $\{M_k = S_k = 2\}_{k=1}^K$ and $N = 4$. We model $\{\mathbf{R}_{bk}, \widetilde{\mathbf{R}}_{mk}\}_{k=1}^K$ as $\{\mathbf{R}_{bk}\}_{k=1}^K = \mathbf{R}_b = \rho_b^{|i-j|}$, $\{\widetilde{\mathbf{R}}_{mk}\}_{k=1}^K = \widetilde{\mathbf{R}}_m = \rho_m^{|i-j|}$, where $0 \leq \rho_b (\rho_m) < 1$. The signal-to-noise-ratio (SNR) is defined as $P_{\max}/\sigma^2$ where $P_{\max}$ is the maximum total BS power and $\sigma^2$ is the noise variance. The SNR is controlled by varying $\sigma^2$ while $P_{\max}$ is set to 10. All simulation results are obtained by averaging over 100 randomly chosen channel realizations.

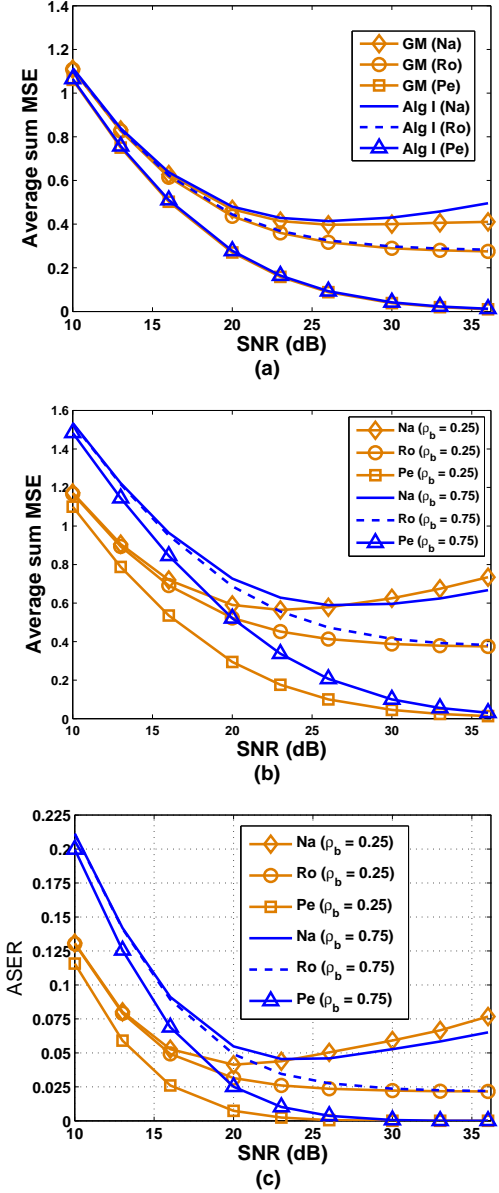

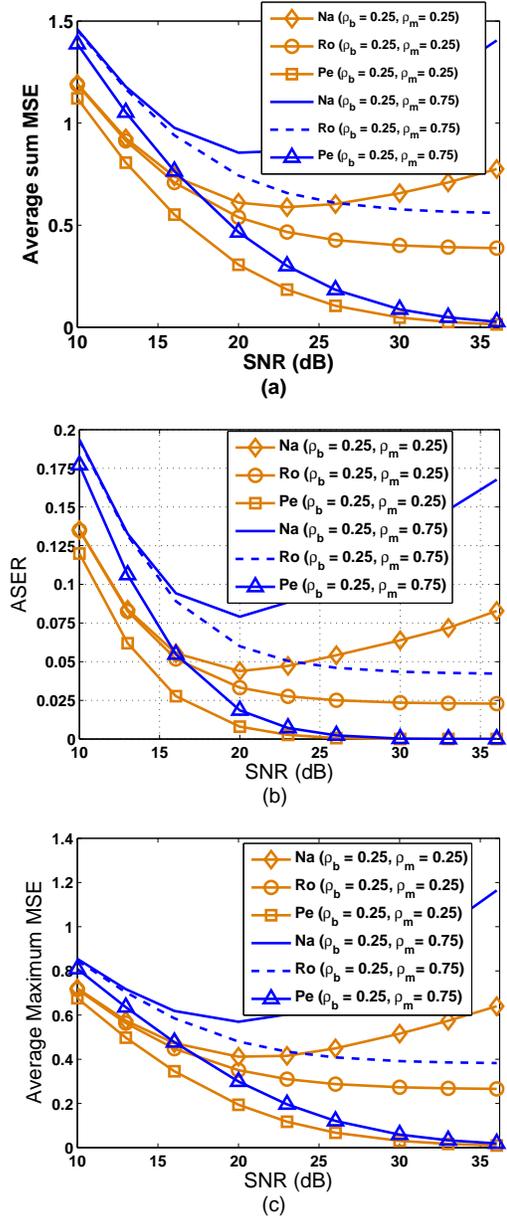

Fig. 2. (a) Comparison of the performance of **Algorithm I** and the GM. (b)-(c) Comparison of the robust and non-robust/naive designs when $\rho_m = 0$, and $\rho_b = 0.25$ and $0.75$. The non-robust/naive, robust and perfect CSI designs are denoted by 'Na', 'Ro' and 'Pe', respectively.

Fig. 3. Comparison of the robust and non-robust/naive designs when $\rho_b = 0.25$, and $\rho_m = 0.25$ and $0.75$. (a)-(b) For the robust sum MSE minimization problem ($\mathcal{P}1$). (c) For the robust Min-max problem ($\mathcal{P}2$).

### A. Simulation results for problem $\mathcal{P}1$

*1) Simulation results for* **Case 1**: For the aforementioned parameters, all our simulation results show that the optimal solution of the SDP problem (the problem (21) after rank relaxation) satisfy the rank constraint of (21)[3]. Consequently, for our setup, the SDP solution is considered as a global minimum (GM) for (17). Similar observation has also been made in the perfect CSI case [4]. Now, we check whether **Algorithm I** achieves the GM or not when $\mathbf{R}_b = \mathbf{I}$, $\sigma_{e1}^2 = \sigma_{e2}^2 = 0.0101$. Fig. 2.a shows that the GM can be achieved by **Algorithm I** for the robust and perfect CSI designs. In the non-robust/naive design, which refers to the design in which the estimated channel is considered as perfect [5], the gap between **Algorithm I** and the GM is large in the high SNR zone.

*2) Comparison of robust and non-robust/naive designs:* For **Case 1**, as can be seen from Fig. 2.a, the robust design has better performance than the non-robust design. Now, we compare the performance of our robust design with that of the non-robust design proposed in [4] for **Case 2**. The comparison is based on the total sum AMSE and the average symbol error rate (ASER)[4] of all users versus the SNR

---

[3]We have noted that the SDP solution of this problem does not always satisfy its rank constraint when $S_k < M_k$. Simulation results for the case $S_k < M_k$ are not included for conciseness.

[4]For the sum AMSE design, ASER is also an appropriate metric for comparing the performance of the robust and non-robust designs [13]. QPSK modulation is utilized for each symbol stream.

when $\{\tau_k = 1\}_{k=1}^K$.

*2.1) When $\{\sigma_{ek}^2 \neq 0\}_{k=1}^K$, $\rho_b \neq 0$ and $\rho_m = 0$:* Here we examine the joint effect of $\rho_b$ and $\{\sigma_{ek}^2\}_{k=1}^K$ on the system performance. To this end, we set $\sigma_{e1}^2 = 0.0101, \sigma_{e2}^2 = 0.0204$ and then we vary $\rho_b$ from 0.25 to 0.75. Figs. 2.(b-c) show that as the BS antenna correlation coefficient increases, the sum AMSE and ASER also increase.

*2.2) When $\{\sigma_{ek}^2 \neq 0\}_{k=1}^K$, $\rho_b \neq 0$, $\rho_m \neq 0$:* Now we discuss the effects of $\{\sigma_{ek}^2\}_{k=1}^K$, $\rho_b$ and $\rho_m$ on the system performance. We keep $\sigma_{e1}^2 = 0.0101, \sigma_{e2}^2 = 0.0204, \rho_b = 0.25$ and then we take $\rho_m$ as 0.25 and 0.75. Figs. 3.(a-b) show that the performance of the system degrades further as $\rho_m$ increases[5].

The results of Section VI-A.2 gracefully fit to that of [13] where (16) is examined for single user MIMO systems.

### B. Simulation results for problem $\mathcal{P}2$

This simulation compares the performance of the robust design and the non-robust design proposed in [19]. Here we keep $\eta_1 = \eta_2 = 0.3$, $\sigma_{e1}^2 = 0.0101, \sigma_{e2}^2 = 0.0204, \rho_b = 0.25$ and then we take $\rho_m$ as 0.25 and 0.75. Fig. 3.c shows that the maximum AMSE of the robust design is less than that of the non-robust design proposed in [19]. Moreover, the performance gap between these designs increases as the SNR increases. This figure also illustrates the fact that large antenna correlation factor degrades the performance of the considered system.

In all figures, the robust design outperforms the non-robust design and the improvement is larger for high SNR regions. This can be seen from the term $\mathbf{\Gamma}_k^{DL}$ of (3) where, at high SNR regions, the effect of $\sigma^2$ is dominated by $\sigma_{ek}^2 \text{tr}\{\mathbf{R}_{bk}\mathbf{G}\mathbf{P}\mathbf{G}^H\}\mathbf{R}_{mk}$ (amplified error). Since the non-robust design does not take into account the effect of $\sigma_{ek}^2 \text{tr}\{\mathbf{R}_{bk}\mathbf{G}\mathbf{P}\mathbf{G}^H\}\mathbf{R}_{mk}$ which is the dominant term, the performance of this design degrades. This implies that as the SNR increases, the performance gap between the robust and non-robust design increases. In all plots, when $\rho_b$ ($\rho_m$) increases, the system performance degrades. This is because as $\rho_b$ ($\rho_m$) increases, the number of symbols with low channel gain increases (this can be easily seen from the eigenvalue decomposition of $\mathbf{R}_{bk}$ ($\mathbf{R}_{mk}$)). Consequently, for fixed total BS power, the total sum AMSE ($\mathcal{P}1$) and maximum AMSE ($\mathcal{P}2$) also increase.

## VII. CONCLUSIONS

In this paper we consider two MSE-based transceiver design problems where imperfect CSI is available both at the BS and MSs. The problems are examined by using the MSE duality approach. The duality are established by transforming only the power allocation matrices from uplink to downlink channel and vice versa. Using our duality results, we propose iterative algorithms that perform optimization alternatively by switching between the uplink and downlink channels. Simulation results show the superior performance of the proposed robust design compared to the non-robust/naive design.

---

[5]**Remark:** When the SNR $\to \infty$ (i.e., $\sigma^2 \to 0$), a sum AMSE floor exists for our robust design. Such sum AMSE floor is observed in Figs. 2.(a-b) and Fig. 3.a. The analytical proof is given as follows: for any $\{\sigma_{ek}^2, \widetilde{\mathbf{R}}_{mk}, \mathbf{R}_{bk}\}_{k=1}^K$ and $\{\tau_k = 1\}_{k=1}^K$, after some mathematical manipulations the uplink sum MAMSE can be expressed as $\text{tr}\{(\mathbf{I}_S + \mathbf{Q}^{1/2}\mathbf{U}^H\widehat{\mathbf{H}}^H(\sigma^2\mathbf{I} + \sum_{i=1}^K \sigma_{ej}^2 \text{tr}\{\mathbf{R}_{mi}\mathbf{U}_i\mathbf{Q}_i\mathbf{U}_i^H\}\mathbf{R}_{bi})^{-1}\widehat{\mathbf{H}}\mathbf{U}\mathbf{Q}^{1/2})^{-1}\}$. Hence, when $\sigma^2 \to 0$ the sum MAMSE approaches to $\text{tr}\{(\mathbf{I}_S + \mathbf{Q}^{1/2}\mathbf{U}^H\widehat{\mathbf{H}}^H(\sum_{i=1}^K \sigma_{ei}^2 \text{tr}\{\mathbf{R}_{mi}\mathbf{U}_i\mathbf{Q}_i\mathbf{U}_i^H\}\mathbf{R}_{bi})^{-1}\widehat{\mathbf{H}}\mathbf{U}\mathbf{Q}^{1/2})^{-1}\} > 0$.